\documentclass[reqno,a4paper,12pt]{amsart}

\usepackage{xcolor}

\usepackage{graphicx}

\usepackage{mathptmx}
\usepackage{mathrsfs}
\usepackage{amsfonts}
\usepackage{amsmath,amsthm,amssymb,amsopn}
\usepackage{amsrefs}
\usepackage{amscd}
\usepackage{cancel}

\usepackage{tikz,pgf,amscd,enumerate,multicol}
\usepackage{amsmath,amssymb,tikz-cd,tkz-graph}%
\usetikzlibrary{decorations.markings,arrows,cd}

\newcommand{\uloopr}[1]{\ar@'{@+{[0,0]+(-4,5)}@+{[0,0]+(0,10)}@+{[0,0] +(4,5)}}^{#1}}
\newcommand{\uloopd}[1]{\ar@'{@+{[0,0]+(5,4)}@+{[0,0]+(10,0)}@+{[0,0]+ (5,-4)}}^{#1}}
\newcommand{\dloopr}[1]{\ar@'{@+{[0,0]+(-4,-5)}@+{[0,0]+(0,-10)}@+{[0, 0]+(4,-5)}}_{#1}}
\newcommand{\dloopd}[1]{\ar@'{@+{[0,0]+(-5,4)}@+{[0,0]+(-10,0)}@+{[0,0 ]+(-5,-4)}}_{#1}}
\newcommand{\luloop}[1]{\ar@'{@+{[0,0]+(-8,2)}@+{[0,0]+(-10,10)}@+{[0, 0]+(2,2)}}^{#1}}

\newcommand{\legendbox}[1]{%
  \textcolor{#1}{\rule{\fontcharht\font`X}{\fontcharht\font`X}}%
}

\usepackage{amsfonts}
\usepackage[mathcal]{eucal}
\usepackage{eufrak}
\usepackage{amssymb}
\usepackage{amsmath}
\usepackage{mathrsfs}
\usepackage[colorlinks]{hyperref}
\usepackage{enumerate}

\definecolor{citecol}{RGB}{145, 1, 1}

\hypersetup{colorlinks=true,citecolor=citecol,linkcolor=blue,linktocpage=true}

\usepackage{tikz}

\usepackage{comment}

\theoremstyle{plain}
\newtheorem {lemma}{Lemma}[section] 
\newtheorem {theorem}[lemma]{Theorem}

\newtheorem {conjecture}{Conjecture}
\newtheorem {question}{Question}

\theoremstyle{definition}
\newtheorem {remark}[lemma]{Remark}

\newtheorem {example}[lemma]{Example}

\theoremstyle{definition}

{}

\newcommand{\M}{\operatorname{\mathbb M}}

\newcommand{\gr}{\operatorname{gr}}

\newcommand{\ES}{\operatorname{ES}}
\newcommand{\Ss}{\operatorname{SS}}

\newcommand{\cO}{\mathcal O}
\def\C{\mathbb{C}} 
\newcommand{\cR}{\mathcal R}
\newcommand{\K}{\mathsf{k}}

\newcommand{\SP}{\operatorname{SP}}

\begin{document}

\title[]{Monoids, dynamics and Leavitt  path algebras}


\author{Gene Abrams}
\address{Department of Mathematics\\University of Colorado\\Colorado Springs, CO 80918, USA}
\email{abrams@math.uccs.edu}

\author{Roozbeh Hazrat}\address{
Centre for Research in Mathematics and Data Science\\
Western Sydney University\\
Australia}
\email{r.hazrat@westernsydney.edu.au}



\begin{abstract}
 Leavitt path algebras, which are algebras associated to  directed graphs, were first introduced about 20 years ago. They have strong connections to such topics as  symbolic dynamics, operator algebras, non-commutative geometry, representation theory, and even chip firing. In this article we invite the reader to   sneak a peek at these fascinating algebras and their interplay with  several seemingly disparate  parts of mathematics. 
\end{abstract}


\maketitle

\tableofcontents

\section{Introduction} \label{introkai}

This article was written for the series {\it Snapshots of modern mathematics from Oberwolfach}.  {\it Leavitt path algebras} were first defined in 2005.    Our intention herein is to give the reader a sense both of how, in the intervening two decades, these structures  have permeated a number of seemingly-unrelated branches of mathematics.   As well, by presenting some open questions in the topic, we hope to entice the reader into further exploring the nature and structure of these algebras.

\section{Chip firing}\label{chipsfire}

We start by looking at a seemingly simple local symmetric behaviour called \emph{chip firing}.  
Although at first glance the chip firing process looks like nothing more than a simple board game, the resulting phenomena can  model extremely complex real world events such as traffic jams or bush fires. 

The mechanics of chip firing are easy to describe.   Suppose at each vertex (node) in a grid there sit some nonnegative numbers of objects ({\it chips}).   Then suppose that  any vertex in the grid transfers ({\it fires})  one  chip to each of its neighbours  if the number of chips on the vertex is at least the degree of the vertex. As an example, in the picture below, we have a grid with four chips sitting on the origin vertex and two on the top middle vertex.  Then the chips on the origin will be fired to the adjacent vertices along the grid. As the picture shows, the top middle vertex now has three chips and has degree three, so it will fire. After this firing, the configuration becomes \emph{stable}, meaning no more firing can occur.

\begin{center}
\begin{tikzpicture}[scale=3.5]
\fill (0,0)  circle[radius=.6pt];
\fill (0,-1/2)  circle[radius=.6pt];
\fill (0,1/2)  circle[radius=.6pt];
\fill (1/2,0)  circle[radius=.6pt];
\fill (-1/2,0)  circle[radius=.6pt];
\fill (-1/2,1/2)  circle[radius=.6pt];
\fill (-1/2,-1/2)  circle[radius=.6pt];
\fill (1/2,1/2)  circle[radius=.6pt];
\fill (1/2,-1/2)  circle[radius=.6pt];
\draw  (0,0) to (0,-1/2);
\draw  (0,0) to (0,1/2);
\draw  (0,0) to (1/2,0);
\draw  (0,0) to (-1/2,0);
\draw  (-1/2,1/2) to (-1/2,0);
\draw  (-1/2,1/2) to (0,1/2);
\draw  (1/2,1/2) to (0,1/2);
\draw  (1/2,1/2) to (1/2,0);
\draw  (-1/2,-1/2) to (0,-1/2);
\draw  (-1/2,-1/2) to (-1/2,0);
\draw  (0,-1/2) to (1/2,-1/2);
\draw  (1/2,-1/2) to (1/2,0);
\draw[->, thick, blue, shorten >=5pt] (0,0) to [out=-40, in=220] (1/2,0);
\draw[->, thick, blue, shorten >=5pt] (0,0) to [out=270 - 40, in=120] (0, -1/2);
\draw[->, thick, blue, shorten >=5pt] (0,0) to [out=120, in=240] (0, 1/2);
\draw[->, thick, blue, shorten >=5pt] (0,0) to [out=210, in=-40] (-1/2, 0);
\draw (0.07,0.07)  node[red, font=\ttfamily\scriptsize] {{\bf 4}};
\draw (0.015,1/2+0.08)  node[red, font=\ttfamily\scriptsize] {{\bf 2}};

\draw[->, double, very thick] (0.62,0) to (0.83, 0);
\draw[->, double, very thick] (1.62+1/2,0) to (1.83+1/2, 0);

\fill (1.5,0)  circle[radius=.6pt];
\fill (1.5,-1/2)  circle[radius=.6pt];
\fill (1.5,1/2)  circle[radius=.6pt];
\fill (1.5+1/2,0)  circle[radius=.6pt];
\fill (1.5-1/2,0)  circle[radius=.6pt];
\fill (1.5 -1/2,1/2)  circle[radius=.6pt];
\fill (1.5-1/2,-1/2)  circle[radius=.6pt];
\fill (1.5+1/2,1/2)  circle[radius=.6pt];
\fill (1.5+1/2,-1/2)  circle[radius=.6pt];
\draw  (1.5+0,0) to (1.5+0,-1/2);
\draw  (1.5+0,0) to (1.5+0,1/2);
\draw  (1.5+0,0) to (1.5+1/2,0);
\draw  (1.5+0,0) to (1.5-1/2,0);
\draw  (1.5-1/2,1/2) to (1.5-1/2,0);
\draw  (1.5-1/2,1/2) to (1.5+0,1/2);
\draw  (1.5+1/2,1/2) to (1.5+0,1/2);
\draw  (1.5+1/2,1/2) to (1.5+1/2,0);
\draw  (1.5-1/2,-1/2) to (1.5+0,-1/2);
\draw  (1.5-1/2,-1/2) to (1.5-1/2,0);
\draw  (1.5+0,-1/2) to (1.5+1/2,-1/2);
\draw  (1.5+1/2,-1/2) to (1.5+1/2,0);
\draw (1.5+0.05+1/2,0.08)  node[red, font=\ttfamily\scriptsize] {{\bf 1}};
\draw (1.5,1/2+0.1)  node[red, font=\ttfamily\scriptsize] {{\bf 3}};
\draw (1.5,-1/2-0.1)  node[red, font=\ttfamily\scriptsize] {{\bf 1}};
\draw (0.9+0.05,0.08)  node[red, font=\ttfamily\scriptsize] {{\bf 1}};
\draw[->, thick, blue, shorten >=5pt] (1.5,1/2) to [out=-40, in=220] (1.5+1/2,1/2);
\draw[->, thick, blue, shorten >=5pt] (1.5,1/2) to [out=220, in=-40] (1.5-1/2, 1/2);
\draw[->, thick, blue, shorten >=5pt] (1.5,1/2) to [out=240, in=120] (1.5, 0);

\fill (1.5+1.5,0)  circle[radius=.6pt];
\fill (1.5+1.5,-1/2)  circle[radius=.6pt];
\fill (1.5+1.5,1/2)  circle[radius=.6pt];
\fill (1.5+1.5+1/2,0)  circle[radius=.6pt];
\fill (1.5+1.5-1/2,0)  circle[radius=.6pt];
\fill (1.5+ 1.5 -1/2,1/2)  circle[radius=.6pt];
\fill (1.5+1.5-1/2,-1/2)  circle[radius=.6pt];
\fill (1.5+1.5+1/2,1/2)  circle[radius=.6pt];
\fill (1.5+1.5+1/2,-1/2)  circle[radius=.6pt];

\draw  (1.5+0+1.5,0) to (1.5+0+1.5,-1/2);
\draw  (1.5+0+1.5,0) to (1.5+0+1.5,1/2);
\draw  (1.5+0+1.5,0) to (1.5+1/2+1.5,0);
\draw  (1.5+0+1.5,0) to (1.5-1/2+1.5,0);
\draw  (1.5-1/2+1.5,1/2) to (1.5-1/2+1.5,0);
\draw  (1.5-1/2+1.5,1/2) to (1.5+0+1.5,1/2);
\draw  (1.5+1/2+1.5,1/2) to (1.5+0+1.5,1/2);
\draw  (1.5+1/2+1.5,1/2) to (1.5+1/2+1.5,0);
\draw  (1.5-1/2+1.5,-1/2) to (1.5+0+1.5,-1/2);
\draw  (1.5-1/2+1.5,-1/2) to (1.5-1/2+1.5,0);
\draw  (1.5+0+1.5,-1/2) to (1.5+1/2+1.5,-1/2);
\draw  (1.5+1/2+1.5,-1/2) to (1.5+1/2+1.5,0);
\draw (1.5+0.08+1/2+1.5,0.08)  node[red, font=\ttfamily\scriptsize] {{\bf 1}};
\draw (3-1/2-0.07,0.08)  node[red, font=\ttfamily\scriptsize] {{\bf 1}};
\draw (1.5+1.5,-1/2-0.1)  node[red, font=\ttfamily\scriptsize] {{\bf 1}};
\draw (3+0.1,0.07)  node[red, font=\ttfamily\scriptsize] {{\bf 1}};
\draw (3+1/2,1/2+0.08)  node[red, font=\ttfamily\scriptsize] {{\bf 1}};
\draw (3-1/2,1/2+0.08)  node[red, font=\ttfamily\scriptsize] {{\bf 1}};

\end{tikzpicture}
\end{center}

The above grid is $3 \times 3$; imagine in its place an analogous $100,000 \times 100,000$ grid.   Now, rather than starting with just 4 chips, suppose we start with $20,000,000$ chips placed  on the origin of the huge grid,  and we start  the firing process. Upon  stabilization, each vertex will have either 0, 1, 2, or 3 chips on it. 
Assigning the colors \legendbox{blue} to the vertices with zero chips,  \legendbox{cyan} with 1 chip,  \legendbox{yellow} for 2 and  \legendbox{brown} for vertices with 3 chips upon stabilisation, it turns out that
 the distribution of the chips form a self-similar symmetric fractal shape (picture below)~\cite{Klivans}. The precise explanation why this fractal shape is formed is yet to be discovered!


\begin{figure}[h]
\includegraphics[width=6.5cm]{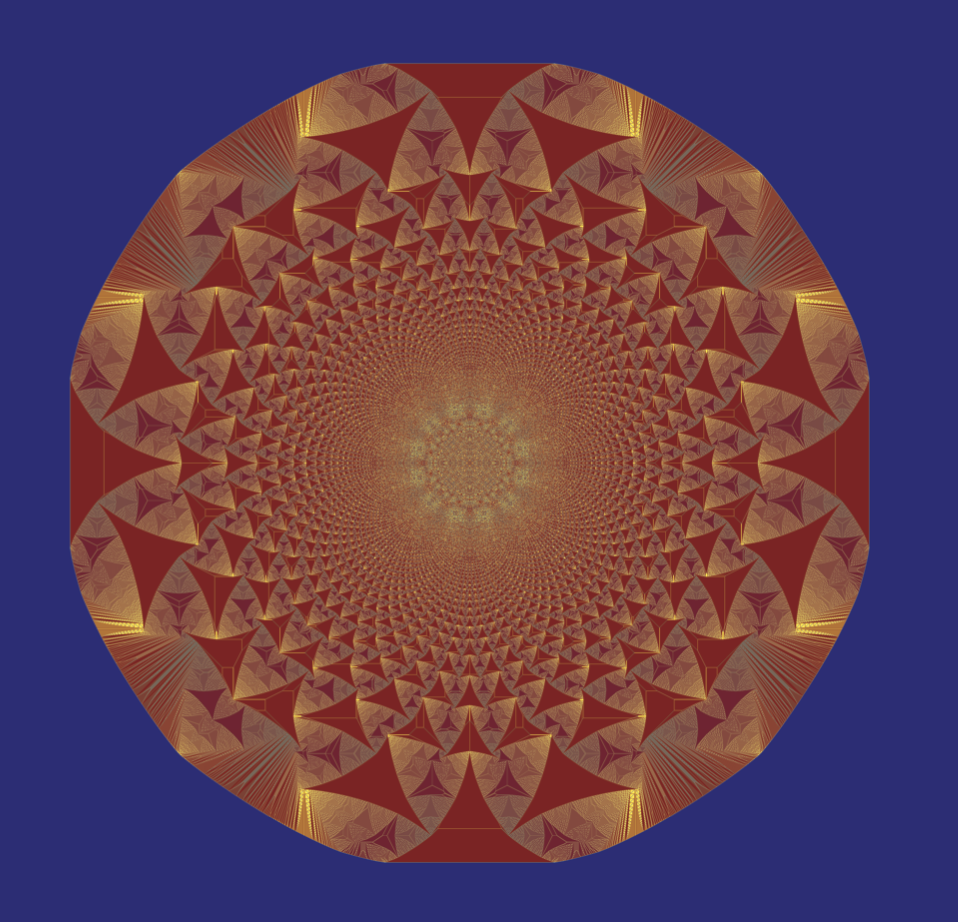}
\end{figure}
The notion of chip firing (also known as a {\it sandpile model}) encapsulates how objects spread and evolve along a grid. The models were conceived in 1987 by Bak, Tang and Wiesenfeld~\cite{bak}  as examples of {\it self-organised criticality}, meaning the tendency of physical systems to organise themselves (without any input from outside the system) toward critical but barely stable states. The models have been used to describe phenomena such as forest fires, traffic jams, stock market fluctuations, etc. The book of Bak~\cite{bakbook} describes how events in nature (strikingly often!)  follow this type of behaviour. 

Next, we enhance the model by using grids that have directions and a ``sink'' that destroys the chips.   
Consider a \emph{sandpile graph}, namely a (finite, directed) graph $E$ with a unique  sink  vertex $s$ such that there is a (directed) path from any vertex of $E$ to $s$ (formal definition in \S\ref{monformal}). Consider a collection of chips placed on each vertex of the sandpile graph (a {\it configuration}). A vertex is {\it unstable} if it has at it the same or more chips than the number of edges emitting from it.  In this case the vertex {\it fires}  by sending one chip along each edge emitting from the vertex to each neighbouring vertex.  This firing may cause neighbouring vertices to become unstable. The assumption is that chips arriving to the vertex $s$ vanish. A configuration is {\it stable} if no vertex is unstable.    Two foundational results used in the theory are:  (1)  the order of firing does not matter
(thus ``abelian'' sandpile models); and  (2) any configuration can be sequentially fired
 to reach a (unique) stable configuration.
 
 Let us make this clearer by presenting the following example.  
 If there are 8 chips on the vertex $v$ of the sandpile graph $E$ below, then the consecutive firings reduce the model into the following distribution.   (The notation is self-evident.)  \label{thispage}
  $$8v \ \leadsto \ 6v+ u  \ \leadsto \ 4v+2u  \ \leadsto \ 2v+3u  \ \leadsto \ 3v+z  \ \leadsto \ v+u+z.$$

\begin{center}
\begin{tikzpicture}[scale=3.5]
\fill (0,0)  circle[radius=.6pt];
\fill (0,-1/2)  circle[radius=.6pt];
\fill (-1/2,-1/2)  circle[radius=.6pt];
\fill (1/2,-1/2)  circle[radius=.6pt];
\draw (0,0.1) node{$u$};

\draw (-1/2+ -0.1,-1/2 + 0.1) node{$v$};
\draw (0,-1/2 + -0.1) node{$s$};
\draw (1/2+ 0.1,-1/2 + 0.1) node{$z$};

\draw (0+0.09,0.1)  node[red, font=\ttfamily\scriptsize] {{\bf 1}};
\draw (-1/2+ -0.1+ 0.15,-1/2 + 0.1)  node[red, font=\ttfamily\scriptsize] {{\bf 1}};
\draw (1/2-0.07,-1/2 + 0.1) node[red, font=\ttfamily\scriptsize] {{\bf 1}};

\draw[->, shorten >=5pt] (0,0) to (0,-1/2);
\draw[->, shorten >=5pt] (-1/2,-1/2) to (0,-1/2);
\draw[->, shorten >=5pt] (1/2,-1/2) to (0,-1/2);
\draw[->, shorten >=5pt]  (1/2,-1/2)to[in=10,out=-70, loop] (1/2,-1/2);
\draw[->, shorten >=5pt] (-1/2,-1/2) to [in=250,out=170, loop] (-1/2,-1/2);
\draw[->, shorten >=5pt] (0,0) to [in=90, out=180] (-1/2,-1/2);
\draw[->, shorten >=5pt] (-1/2,-1/2) to [in=250, out=10] (0,0);
\draw[->, shorten >=5pt] (0,0) to [in=90, out=0] (1/2,-1/2);
\draw[->, shorten >=5pt] (1/2,-1/2) to [in=290, out=170] (0,0);

\draw[->, double, very thick] (-1.1,-1/4) to (-1.1+0.21, -1/4);

\fill (-2+0,0)  circle[radius=.6pt];
\fill (-2+0,-1/2)  circle[radius=.6pt];
\fill (-2 -1/2,-1/2)  circle[radius=.6pt];
\fill (-2 + 1/2,-1/2)  circle[radius=.6pt];
\draw (-2 ,0.1) node{$u$};
\draw (-2-0.8,-0.05) node{$E:$};

\draw (-2-1/2+ -0.1,-1/2 + 0.1) node{$v$};
\draw (-2-1/2+ 0.05,-1/2 + 0.1) node[red, font=\ttfamily\scriptsize] {{\bf 8}};

\draw (-2,-1/2 + -0.1) node{$s$};
\draw (-2+1/2+ 0.1,-1/2 + 0.1) node{$z$};
\draw[->, shorten >=5pt] (-2,0) to (-2,-1/2);
\draw[->, shorten >=5pt] (-2-1/2,-1/2) to (-2,-1/2);
\draw[->, shorten >=5pt] (-2+1/2,-1/2) to (-2+0,-1/2);
\draw[->, shorten >=5pt]  (-2+1/2,-1/2)to[in=10,out=-70, loop] (-2+1/2,-1/2);
\draw[->, shorten >=5pt] (-2-1/2,-1/2) to [in=250,out=170, loop] (-2-1/2,-1/2);
\draw[->, shorten >=5pt] (-2,0) to [in=90, out=180] (-2-1/2,-1/2);
\draw[->, shorten >=5pt] (-2-1/2,-1/2) to [in=250, out=10] (-2,0);
\draw[->, shorten >=5pt] (-2,0) to [in=90, out=0] (-2+1/2,-1/2);
\draw[->, shorten >=5pt] (-2+1/2,-1/2) to [in=290, out=170] (-2,0);

\end{tikzpicture}
\end{center}
Many questions arise, e.g.: How many firings are needed for a configuration to be stabilised?

The theory of chip firing has found many applications, and the  books  \cites{Klivans, corper}   are dedicated to this theory.  We return to this model after introducing two mathematical objects:  {\it monoids}  and {\it Leavitt path algebras}.

\section{Symbolic dynamics}


According to \cite{lindmarcus}, the subject of symbolic dynamics began in 1898 when Jacques Hadamard in \cite{Had}  had the idea  of representing dynamical systems (e.g., motion of planets, or, in Hadamard's case,  geodesic flows on surfaces of negative curvature) by discrete sequences of information. The main point of Hadamard's work is that there is often a relatively simple description of the possible sequences that can arise in such a representation.   More generally,   symbolic dynamics deals with sequences (often indexed by $\mathbb{Z}$) of symbols chosen from a finite set,  equipped with a \emph{shift map} which creates the dynamical behaviour.    One of the fundamental objects of study in symbolic dynamics is called a \emph{shift of finite type} (SFT), which consists of those sequences that do not include certain ``forbidden'' finite sequences of symbols taken from the underlying finite set. The applications  abound, from topological quantum field theory, ergodic
theory, and statistical mechanics to coding and information theory~\cites{Kim-RoushW1999,lindmarcus}. 

The fundamental work of Williams~\cite{williams} shows that understanding these dynamical systems reduces to understanding certain relations between matrices. (Echoing T. Y. Lam \cite[p.11]{Lam} here: ``God bless matrices!")  
Thus, instead of talking about ``conjugacy'' between shifts of finite type (a precise topological notion), we can talk about ``strong shift equivalence'' between matrices. We provide the definitions here.   (Additional information can be found in \cite{lindmarcus}.)

Two square nonnegative integer matrices $A$ and $B$ are called {\it elementary shift equivalent},  denoted  $A\sim_{\ES} B$, if there are nonnegative integer matrices $R$ and $S$ such that 
$$A=RS\,\, \, \text{and } \, \, B=SR.\\$$
Note that $R$ and $S$ need not be square matrices.

\begin{example}
The square matrices $A=(2)$ and 
$B=\begin{pmatrix}
1 & 1 \\
1& 1\\
\end{pmatrix}$ are elementary shift equivalent, because 
$$
(2)=
\begin{pmatrix}
1 & 1 \\
\end{pmatrix} 
\begin{pmatrix}
1  \\
1\\
\end{pmatrix},  \ \ \ \mbox{and}  \ \ \ 
\begin{pmatrix}
1 & 1 \\
1& 1\\
\end{pmatrix}
=
\begin{pmatrix}
1  \\
1\\
\end{pmatrix}
\begin{pmatrix}
1 & 1 \\
\end{pmatrix}.$$

\end{example}

The terminology ``elementary shift equivalent" is standard, but somewhat misleading:  $\sim_{\ES}$ is {\it not} an equivalence relation (it is not transitive).   The equivalence relation $\sim_{\Ss}$  generated by elementary shift equivalence on the set of square nonnegative integer matrices (of any finite size) is called {\it strong shift equivalence}. That is,  $A\sim_{\Ss} B$ in case there is a sequence of square matrices $A_i$ such that $$A\sim_{\ES} A_1 \sim_{\ES} A_2 \sim_{\ES}  \cdots  \sim_{\ES} A_k \sim_{\ES} B.$$

It is extremely hard in practice to determine whether  two matrices are strong shift equivalent (equivalently, whether two SFT's are conjugate).   As an example, the matrices 
$A=\begin{pmatrix}
1 & 3 \\
2& 1\\
\end{pmatrix}$ and $B=\begin{pmatrix}
1 & 6 \\
1& 1\\
\end{pmatrix}$ are strong shift equivalent (a fact that was established using the brute force of computer computations)~\cite{lindmarcus}.   However, the following seemingly-basic question is yet to be answered. 

 \begin{question}

Consider the matrices $A_k=\begin{pmatrix}
1 & k \\
k-1 &1\\
\end{pmatrix}$ and $B_k=\begin{pmatrix}
1 & k(k-1)\\
1 & 1\\
\end{pmatrix},$ \  for integers $k\geq 3$. As indicated directly above, after some very heavy lifting  it was  established  that $A_3\sim_{\Ss} B_3.$  
Is  $A_k \sim_{\Ss} B_k$ for $k\geq 4$?

\end{question}


Next we introduce the notion of shift equivalence.  The nonnegative  integer square matrices $A$ and $B$ are called {\it shift equivalent}, denoted $A\sim_S B$,  if there exist $\ell\ge 1$ and nonnegative integer matrices $R$ and $S$ such that 
\begin{gather}
A^\ell= RS, \ \ \, B^\ell=SR, \label{eq:equidentro}\\
AR=RB, \ \ \, SA=BS.\label{eq:equifuera}
\end{gather}

When $\ell=1$,  it is easy to show (just using associativity of multiplication of matrices) that  equations (1) imply equations  (2), so the $\ell = 1$ case of shift equivalence is the same as elementary shift equivalence.    By an easy induction it follows 
that strong shift equivalence implies shift equivalence. 

The following ``theorem" was published in 1974  in {\it Annals of Mathematics}, one of the most prestigious mathematics journals.      

\begin{theorem}[Williams] \cite{williams}   \label{willwrong}
Let $A$ and $B$ be nonnegative square integral matrices. Then $A\sim_{S} B$ \xcancel{if and only if}   $A\sim_{\Ss} B$.
\end{theorem}

It turned out that Williams' proof that shift equivalence implies strong shift equivalence 
was not correct. 
Once the error in the ``proof" of that implication was realised, the statement remained as a ``conjecture'' for about 20 years. But it was finally {\it disproved}, via  counterexamples presented by Kim and Roush  in \cites{kimroush}. The extremely subtle  distinction between shift equivalence and strong shift equivalence remains elusive. Indeed, the following two questions remain active areas of current research:  (1) to determine the decidability of strong shift equivalence (i.e., to obtain an algorithm that determines whether two matrices are strong shift equivalent, or to show that no such algorithm exists);  and (2)  to establish useful sufficient  conditions that yield that   strong shift equivalence and shift equivalence coincide.    It is interesting to note that the property ``shift equivalent" has been shown to be decidable.  

We return to the notions of shift equivalence and strong shift equivalence after introducing  {\it monoids} and {\it Leavitt path algebras}.

\section{Commutative Monoids} 

Perhaps the most ubiquitous  of all algebraic structures are commutative monoids. A nonempty set $M$ with a binary operation $+$ is called a \emph{commutative monoid} if for all $a,b,c\in M$ we have $a+b=b+a$, $a+(b+c)= (a+b)+c$, and there exists $0\in M$ such that $a+0=a$. The first two conditions guarantee that we obtain the same result irrelevant of whether reading the expressions from left to right or from right to left. Examples of commutative monoids are abundant:  for instance,  let $X$ be a nonempty set and $\mathcal P(X)$ its power set;  with  either union or intersection  as a binary operation, $\mathcal P(X)$ is a commutative monoid. 

For a unital ring $A$, the isomorphism classes of finitely generated projective left (or right)  $A$-modules with direct sum as the operation form a commutative monoid,  denoted by $\mathcal V(A)$. 
(This construction can naturally be extended to non-unital rings via idempotents.) Explicitly, let $[P]$ denote the class of $A$-modules isomorphic to $P$. Then   
the set 
\begin{equation}\label{zhongshan}
\mathcal V(A)=\big \{ \, [P] \mid  P  \text{ is a finitely generated projective A-module} \, \big \}
\end{equation}
 with addition $[P]+[Q]=[P\bigoplus Q]$ forms a commutative monoid.

For an abelian group $\Gamma$ and $\Gamma$-graded unital ring $A$, considering the graded finitely generated projective modules, it provides us with the commutative monoid  $\mathcal V^{\gr}(A)$ which has an action of $\Gamma$ on it via the shift operation on modules (see \cite[Section 3]{hazbook} for the general theory).


\subsection{Graph monoids, sandpile monoids, and talented monoids} \label{monformal}

There are two very natural commutative monoids that can be associated to any directed graph $E$:  the \emph{graph monoid} $M_E$,  and the \emph{talented monoid} $T_E$. We will see that these monoids provide a bridge linking all the previous topics to the theory of Leavitt path algebras. 

Let $E = (E^0, E^1, r,s)$ be a \emph{directed graph}. Here $E^{0}$ and $E^{1}$ are
sets and $r,s$ are maps from \emph{edge set} $E^1$ to \emph{vertex set} $E^0$.  We think of each $e \in E^1$ as an edge 
pointing from $s(e)$ to $r(e)$. A vertex $v\in E^0$ is a sink if $s^{-1}(v)=\emptyset$. 
For example for the graph $E$ on page \pageref{thispage}, $E^0=\{u,v,z,s\}$ with $s$ being the only sink vertex,  while $E^1$ consists of the indicated nine edges.
A \emph{path} $p$ in $E$ is
a sequence $p=e_{1}e_{2}\cdots e_{n}$ of edges $e_{i}$ in $E$ such that
$r(e_{i})=s(e_{i+1})$ for $1\leq i\leq n-1$. For the path $p$, we write $s(p)=s(e_1)$ and $r(p)=r(e_n)$.

Let $w:E^0\rightarrow \mathbb N = \{0,1,2, \dots \}$ be a {\it weight} map, which assigns a non-negative integer to each vertex.  (One example of such a weight map is the {\it vertex weighting}, in which $w(v) = | s^{-1}(v)|$, i.e., $w(v)$ is  the number of edges having $v$ as their source vertex.) 

We define the \emph{graph monoid} to be the commutative monoid

\begin{equation}\label{graphmon}
M_{(E,w)}= \Big \langle \,  v \in E^0 \, \,  \Big \vert \, \, \ w(v) v= \sum_{e\in s^{-1}(v)} r(e) \, \Big \rangle,
\end{equation}
where the relations are taken over vertices $v$ which are not sinks. If the weight of every vertex in $E^0$ is $1$, the graph monoid of (\ref{graphmon}) reduces to 

\begin{equation}\label{graphmonm}
M_E= \Big \langle \,  v \in E^0 \, \,  \Big \vert \, \, v= \sum_{e\in s^{-1}(v)} r(e) \, \Big \rangle . \\
\end{equation}
\noindent
To explicitly describe these graph monoids, we generate the free monoid with basis $E^0$, i.e., the direct sum $\bigoplus_{u\in E^0} N_u$, where $N_u= \mathbb N$,  subject to the congruence relation generated by the relations $w(v)\mathbf e_v=\sum_{e\in s^{-1}(v)}\mathbf e_{r(e)}$. Here $\mathbf e_v \in \bigoplus_{u\in E^0} N_u$ is a row vector  with $1$ in the $v$ entry and zero elsewhere.  

Sandpile monoids, which model the process of chip firing described in Section~\ref{chipsfire}, are  variations of graph monoids.  A finite directed graph $E$ is called a \emph{sandpile graph} if $E$ has a unique sink (denote it by $s$), and for every $v \neq s \in E^0$ there is a  path $p$ with $s(p)=v$ and $r(p) = s$.  
 Then the \emph{sandpile monoid} $\SP(E)$ of the sandpile graph $E$  is the monoid defined by generators and relations as
\begin{equation}\label{sandexpi}
\SP(E)  =   \Big \langle   v \in E^0 \, \,  \Big \vert \, \,  s = 0;  \     |s^{-1}(v)| v= \sum_{e \in s^{-1}(v)} r(e), \ \mbox{for} \  v\in E^0\ \setminus\{s\} \Big \rangle.
\end{equation}
\noindent
Recast, $\SP(E)$ is the monoid $M_{(E,w)}$ (where $w$ is the vertex-weighting), with the additional relation that $s=0$.

Finally, for any graph $E$ we define the \emph{talented monoid} $T_E$, as follows: 

\begin{equation}\label{talentedmon}
T_E= \Big \langle \, v(i), v \in E^0, i \in \mathbb Z  \, \,  \Big \vert \, \, v(i)= \sum_{e\in s^{-1}(v)} r(e)(i+1) \, \Big \rangle,
\end{equation}
where again the relations are taken over vertices $v$ which are not sinks.  
The monoid $T_E$ is equipped with a $\mathbb Z$-action: $n\in \mathbb Z$ acts on the generators by $${}^n v(i):=v(n+i),$$ and the action is extended to all elements of $T_E$,  $\mathbb Z$-linearly. We say that talented monoids $T_E$ and $T_F$ are \emph{$\mathbb Z$-isomorphic}, if there is a monoid isomorphism $\phi:T_E\rightarrow T_F$, which also preserves the $\mathbb Z$-action, i.e., $\phi({}^n a) = {}^n \phi(a)$, for all $a\in T_E$ and $n\in \mathbb Z$.

The relations in $M_E$ can be thought of  as saying that the ``worth'' of the vertex $v$ is the sum of the worths of  the adjacent vertices that $v$ is connected to (counting the weights of adjacent vertices multiple times if there are multiple edges from $v$ to a vertex). Another interpretation is that $v$ ``transforms'' to the sum of the adjacent vertices. With this interpretation, the relations in $T_E$ incorporate a time evolution in the transformations. The following example shows, despite the similarity of the definitions of the monoids $M_E$ and $T_E$, they are in general quite different. 

\begin{example}

Consider the graph $E$ 
\begin{center}
 \begin{tikzpicture}[scale=3.5]
\fill (0,0)  circle[radius=.6pt];
\draw (0,-.2) node{$u$};
\fill (.5,0)  circle[radius=.6pt];
\draw (.5,-.2) node{$v$};
\draw[->, shorten >=5pt] (0,0) to[in=135, out=-135, loop] (0,0);
\draw[->, shorten >=5pt] (0,0) to[in=120,out=60] (.5,0);
\draw[->, shorten >=5pt] (.5,0) to[in=-60,out=-120] (0,0);
\end{tikzpicture}
\end{center}
whose adjacency matrix is $A_E= \left(\begin{matrix} 1 & 1\\ 1 & 0
\end{matrix}\right)$.  
One can  show that $T_E$ is the direct limit of the chain of monoids 
\begin{equation} \label{talmonse}
\mathbb N\oplus \mathbb N
\stackrel{A_E}{\longrightarrow}  \mathbb N\oplus \mathbb N
\stackrel{A_E}\longrightarrow  \mathbb N\oplus \mathbb N
\stackrel{A_E}\longrightarrow \cdots.
\end{equation}

\noindent
This monoid is explicitly described  in \cite[Example IV.3.6]{davidson} as
\[T_E \ = \ \varinjlim_{A_E} \{ (\mathbb N \oplus \mathbb N, A_E) \} \ = \ \Big \{ (m,n)\in \mathbb Z \oplus \mathbb Z  \, \, \Big \vert \, \,   \frac{1+\sqrt{5}}{2} m +n \geq 0 \Big \},\]
  with the action 
\[{}^1 (m,n)=(m,n) \left(\begin{matrix} 1 & 1\\ 1 & 0
\end{matrix}\right)= (m+n,m).\]

\noindent
In contrast, $M_E \cong \{0, x\}$, where $x=x+x$. 

If we modify $E$ by adding  an edge from the vertex $v$ to a sink $s$, the result is a sandpile graph $F$, and the sandpile monoid of (\ref{sandexpi}) becomes $\SP(F) \cong \{0,x,2x,3x\}$ with $3x=4x$. 

\end{example}

\begin{example}

We give one more example to demonstrate the richness of the talented monoid. Consider the following graph $E$: 
\begin{center}
\begin{tikzpicture}[scale=3.5]
\fill (0,0)  circle[radius=.6pt];
\fill (-.35,-.61)  circle[radius=.6pt];
\fill (.35,-.61)  circle[radius=.6pt];
\draw[->, shorten >=5pt] (0,0) to (.35,-.61);
\draw[->, shorten >=5pt] (.35,-.61) to (-.35,-.61);
\draw[->, shorten >=5pt] (-.35,-.61) to [in=-120,out=-60] (.35,-.61) ;
\draw[->, shorten >=5pt] (-.35,-.61) to (0,0);
\draw[->, shorten >=5pt] (0,0) to[in=130,out=50, loop] (0,0);
\draw[->, shorten >=5pt] (-.35,-.61) to[in=250,out=170, loop] (-.35,-.61);
\draw[->, shorten >=5pt] (.35,-.61) to[in=10,out=-70, loop] (-.35,-.61);
\end{tikzpicture}
\end{center}
It is easy to calculate that $M_E \cong \{0,x\}$ with $x+x=x$.   
 However, using deep linear algebra results one can show that 
$T_E=\left\{ {\bf u} \in \mathbb{Z}^3 \, | \, {\bf u} \cdot {\bf z} \geq 0 \right\},$ where $\bf{z}$ is the column vector 
$$\left[
\begin{array}{c}
 \frac{1}{3} \left(1-5 \sqrt[3]{\frac{2}{11+3
   \sqrt{69}}}+\sqrt[3]{\frac{1}{2} \left(11+3
   \sqrt{69}\right)}\right) \\ 
   \frac{1}{3}
   \left(-1+\sqrt[3]{\frac{25}{2}-\frac{3
   \sqrt{69}}{2}}+\sqrt[3]{\frac{1}{2} \left(25+3
   \sqrt{69}\right)}\right) \\
   1 \\
\end{array}
\right],$$ and where  
$${}^1(a,b,c)=(a+c,a+b+c, b+c).$$

\end{example}

For completeness of this section we include the following short discussion (which admittedly  ramps up the notational complexity of our article significantly, albeit temporarily).  
Similar to (\ref{talmonse}),  a nonnegative integral $n\times n$ matrix $A$ gives rise to a direct system of free abelian groups with $A$ acting as an order preserving group homomorphism
\[\mathbb Z^n \stackrel{A}{\longrightarrow} \mathbb Z^n \stackrel{A}{\longrightarrow}  \mathbb Z^n \stackrel{A}{\longrightarrow} \cdots . 
\] 
We regard $\mathbb Z^n$ as a partially ordered group with positive cone $\mathbb N^n$. The direct limit of this system, $\Delta_A:= \varinjlim_{A} \mathbb Z^n$  is a partially ordered group whose positive cone $\Delta^+$ is the direct limit of the associated direct system of positive cones. Multiplication by $A$ induces an automorphism $\delta_A$ of partially ordered groups. The triple  $(\Delta_A, \Delta_A^+, \delta_A)$ is called \emph{Krieger's dimension group}. \index{Krieger's dimension group}

The following theorem was proved by Krieger (\cite{krieger}*{Theorem~4.2}; see also~\cite{lindmarcus}*{Section 7.5} for a detailed algebraic treatment). 

\begin{theorem}[Krieger]
\label{kriegerthm}
Let $A$ and $B$ be  square nonnegative integer matrices. Then the following are equivalent:

 \begin{enumerate}[\upshape(1)]
\item $A$ and $B$ are shift equivalent;

\item there is an isomorphism $(\Delta_A, \Delta_A^+, \delta_A) \cong (\Delta_B, \Delta_B^+, \delta_B);$

\item there is a $\mathbb Z$-isomorphism of talented monoids $T_A \cong T_B$.

\end{enumerate}

\end{theorem}

Here $T_A$ is the talented monoid of the graph whose adjacency matrix is $A$.   
Theorem~\ref{kriegerthm} shows that both Krieger's dimension groups, as well as  talented monoids, are complete invariants for the notion of shift equivalence. It is yet an unsolved problem whether one can similarly  find a complete invariant for strong shift equivalence. In the next section, we will see that Leavitt path algebras can play a role in this direction.

\section{Algebras without unique dimension: the Leavitt algebras $L_\K(m,n)$} 

The supposition that \emph{spaces} have a well-defined (i.e., unique) {\it dimension} is crucial for many results that we expect and use in mathematics and physics. For instance, the fact that if two distinct lines in the standard real plane intersect, they intersect at one point  (or, equivalently, if a system of linear equations with coefficients in $\mathbb{R}$ in two variables has solutions,  the solution is unique)  follows from the  fact that real vector spaces do have unique dimensions. 
The uniqueness of dimension translates to saying that  if two  vector spaces defined over a field $\K$, of dimensions $n$ and $m$, are isomorphic (i.e, if $\K^n\cong \K^m$), then $n=m$. Another interpretation of this property of fields is that only square matrices can be  invertible: namely, if $A\in \M_{m\times n}(\K)$ and $B\in \M_{n \times m}(\K)$ such that $AB=I_m$ and $BA=I_n$ then $m=n$. (Here 
$I_\ell$ is the $\ell \times \ell$ identity matrix, in which all entries are zero, except on diagonals, which are 1.) 

The fact that nonzero elements in a field are invertible can be used to show that vector spaces (over fields)  have a unique dimension. One can  show that the dimensions over other more general rings are also unique: for example, if there is a unital ring homomorphism  $\phi: R\rightarrow \K$, then  spaces over $R$ have a unique dimension.  \emph{Proof}: considering invertible rectangular matrices over $R$, their images in $\K$ are also invertible, so they must be square.   In particular, dimensions over $\mathbb{Z}$ are unique.  

In the 1950's, William Leavitt of the University of Nebraska studied, in a universal manner, rings which do not have unique dimensions. Namely, there are invertible rectangular matrices with entries from these rings which are not square matrices. Leavitt's construction shows precisely how this behaviour occurs.   Let $\K$ be a field  and let $A=(x_{ij})$ and $B=(y_{ji})$, where $1\leq i \leq m$ and $1\leq j \leq n$, be a collection of symbols considered as $m\times n$ and $n\times m$ matrices, respectively.  Then Leavitt's algebras are defined as follows:  
\begin{equation}\label{firstdef}
L_\K(m,n) \  := \  \K\big\langle A, B\big\rangle\big /  \big\langle A B = I_m,  \ BA=I_n \big\rangle.
\end{equation}
\noindent
Here $\K\langle A, B\rangle$ is the free associative (noncommutative) unital ring generated by symbols $x_{ij}$ and $y_{ji}$'s with coefficients from the field $\K$. Furthermore,  $AB=I_m$ and $BA=I_n$ stand for the collection of relations after multiplying the matrices and comparing the two sides of the equations. 

Let us look at the case where $m=1$, which Leavitt first studied in his paper \cite{vitt62}. The algebra (\ref{firstdef}) is  then 
 the free unital associative $\K$-algebra $L_n=L_\K(1,n)$ generated by symbols $\{x_i,y_i \mid 1\leq i \leq n\}$ subject to the following relations: 
\begin{equation}\label{jh54320}
y_ix_j =\delta_{ij}1_\K  \  \text{ for all } 1\leq i,j \leq n, \   \text{  and  }  \ \sum_{i=1}^n x_iy_i=1_\K.
\end{equation} 
(Here $\delta_{ij}$ is Kronecker delta.)  The relations guarantee that the right $L_n$-module homomorphism 
\begin{align}\label{is329ho}
\phi:L_n&\longrightarrow L_n^n\\
a &\mapsto (y_1a	,y_2a,\dots,y_na)\notag
\end{align}
has an inverse 
\begin{align}\label{is329ho9}
\psi:L_n^n&\longrightarrow L_n\\
(a_1,\dots,a_n) &\mapsto  x_1a_1+\dots+x_na_n, \notag 
\end{align}
so $L_n\cong L_n^n$ as right $L_n$-modules. Leavitt showed that 
$n$ is the smallest $k\ge 2$ such that $L_n\cong L_n^k$, and that $L_n$ is a simple ring. In fact, Leavitt proved that these rings are ``super simple":   for any nonzero $x\in L_n$ there exist $a,b\in L_n$ such that
\begin{equation}\label{eq:lnspi}
axb=1.    
\end{equation}
(The point here being that one need not use sums of expressions involving the nonzero element $x$ in order to generate the identity element of the ring.)  Note, however, that $L_n$ is not a division ring or even a domain; in fact we have 
$y_2x_1=0.$     

It turns out that the relations (\ref{jh54320}) are quite ubiquitous in mathematics, and as a result  these Leavitt algebras have fascinating properties and appear in many areas. 
For example Cuntz' algebra $\cO_n$, introduced in \cite{On}, is a $C^*$-algebra completion of $L_\C(1,n)$. 


The ring $L_n$ is already very interesting even when $n=2$.  We have 
\begin{equation}\label{free2free}
L_2 \cong L_2 \oplus L_2
\end{equation}
as $L_2$-modules.   Additionally,  for any $\K$-algebra $A$ having  $A \cong A\oplus A$ as $A$-modules,  it is not hard to show that there is a ring monomorphism $L_2 \rightarrow A$.  Many fundamental questions about this easy-to-describe algebra are yet to be answered.

\begin{question}\label{L2mother} 
Is $L_2$ the mother of all rings?  More formally:  
 do all $\K$-algebras with countable basis embed in $L_2$?  
 (See Remark \ref{L2motherRemark} below.)  
\end{question}

\begin{question}
 Can one explicitly describe all irreducible (resp., indecomposable) modules over  $L_2$?  
  \end{question}

\section{Algebras from graphs:  Leavitt path algebras, and beyond}

In the landmark work \cite{kum}, Kumjian, Pask and Raeburn introduced the $C^*$-algebra $C^*(E)$ for a row-finite graph $E$.   (This work was a natural next step in the development of  ideas presented by Cuntz and Krieger in 
 \cite{ck}.)   
In this context, the Cuntz algebra $\cO_n$ mentioned above is realized as $C^*(\cR_n)$, where 
$\cR_n$ is the  ``rose with $n$-petals" graph
\begin{center}
\begin{tikzpicture}[scale=3.5]
\fill (0,0)  circle[radius=.6pt];
\draw (0.35,-0.1) node{$x_1$};
\draw (0.3,0.2) node{$x_2$};
\draw (-0.09,0.33) node{$x_3$};
\draw (-0.09,-0.33) node{$x_n$};

\draw[->, shorten >=5pt] (0,0) to[out=0, in=90, loop] (0,0);
\draw[->, shorten >=5pt] (0,0) to[out=0+50, in=90+50, loop] (0,0);
\draw[->, dotted, shorten >=5pt] (0,0) to[out=0+2*60, in=90+2*60, loop] (0,0);
\draw[->,  shorten >=5pt] (0,0) to[out=0+3*60, in=90+3*60, loop] (0,0);
\draw[->, shorten >=5pt] (0,0) to[out=0+5*60, in=90+5*60, loop] (0,0);

\end{tikzpicture}
\end{center}

The construction was later further generalised to all graphs in \cite{dritom}.  
The Leavitt path algebra $L_\K(E)$ over a field $\K$ was introduced in \cite{aap05} and \cite{amp} as the direct algebraic counterpart to $C^*(E)$.  (See \cite{TheBook} for the details of the construction of $L_\K(E)$.)  Leavitt's algebras are  recovered in this more general context as $L_\K(1,n) \cong L_\K(\cR_n)$ for each $n>1$. 

The definition of these algebras is as follows.  Let $E = (E^0,E^1, s, r)$ be a directed graph and $\K$ a field. The \emph{Leavitt path algebra} $L_{\K}(E)$ is the universal $\K$-algebra generated by formal elements $\{ v \ | \ v \in E^0 \}$ and $\{  e,e^* \ | \ e\in E^1 \ \}$ subject to the relations
\begin{itemize}
    \item $v^2 = v$, and $vw = 0$ if $v\neq w$ are in $E^0$;
    \item $s(e) e = e r(e) = e, \ r(e)e^* = e^*s(e) = e^*$, for all $e\in E^1$;
    \item $e^*e = r(e)$, for all $e\in E^1$;
    \item $v= \sum_{e\in s^{-1}(v)}ee^*$, for each $v\in E^0$ which is  neither a  sink nor an infinite emitter.
\end{itemize}

Abrams and Aranda Pino in~\cite{aap05} coined the name ``Leavitt path algebra'' and gave the first theorem in the theory.   Concentrating here on finite graphs, the theorem takes the following form.  (A {\it cycle} in a graph is a path which starts and ends at the same vertex, and for which no vertex is repeated.)

\begin{theorem}\label{firsttheorem}
Let $E$ be a finite graph and $\K$ any field. The Leavitt path algebra $L_\K(E)$ is a simple ring if and only if 
\begin{enumerate}[\upshape(1)]
\item  Every vertex in $E$ connects to every cycle and every sink of $E$, and 
\item Every cycle in $E$ has an exit.   (This means that for each cycle $c=e_1e_2\cdots e_n$ in $E$ there exists an edge $e_i$ and an edge $f\neq e_i$ in $E$ for which $s(e_i) = s(f)$.)
\end{enumerate}
\end{theorem}

\begin{example}
Among the graphs below only the one on the right satisfies all the conditions of Theorem~\ref{firsttheorem} and thus gives rise to a simple Leavitt path algebra. 

\begin{center}
\begin{tikzpicture}[scale=3.5]
\fill (-0.5,0)  circle[radius=.6pt];
\fill (0,0)  circle[radius=.6pt];
\draw[->, shorten >=5pt] (-0.5,0) to[in=135, out=-135, loop] (0,0);
\draw[->, shorten >=5pt] (-0.5,0) to[in=120,out=60] (0,0);
\draw[->, shorten >=5pt] (0,0) to[in=-60,out=-120] (-0.5,0);

\fill (-3,0)  circle[radius=.6pt];
\fill (-2.5,0)  circle[radius=.6pt];
\fill (-3.35,0)  circle[radius=.6pt];
\draw[->, shorten >=5pt] (-3.35,0) to (-3,0);
\draw[->, shorten >=5pt] (-3,0) to[in=120,out=60] (-2.5,0);
\draw[->, shorten >=5pt] (-2.5,0) to[in=-60,out=-120] (-3,0);

\fill (-2,0)  circle[radius=.6pt];
\fill (-1.5,0)  circle[radius=.6pt];
\fill (-1.15,0)  circle[radius=.6pt];

\draw[->, shorten >=5pt] (-1.5,0) to (-1.15,0);
\draw[->, shorten >=5pt] (-2,0) to[in=120,out=60] (-1.5,0);
\draw[->, shorten >=5pt] (-1.5,0) to[in=-60,out=-120] (-2,0);
\end{tikzpicture}
\end{center}
\end{example}

\smallskip

Since the publication of \cite{aap05} there has been an explosion of activity in pursuit of  establishing  correspondences between ring-theoretic properties of $L_\K(E)$  and combinatorial properties of $E$.  Many of those correspondences led to new examples and counterexamples in ring theory. Similar explorations are yet to be completed for the relationship between monoid properties of $T_E$ and combinatorial properties of $E$. We give here one sample of a situation in which the three  structures $T_E$, $E$, and $L_\K(E)$  are clearly tightly related. Recall that a ring $R$ is called \emph{Zorn} if every  element $a\in R$ is either nilpotent (i.e., $a^n=0$ for some $n \in \mathbb{N}$) or there exists $b\in R$ such that $ab$ is a nonzero idempotent. The following theorem is a combination of results from~\cite{rangazorn} and \cite{hazli}.

\begin{theorem}
Let $E$ be a finite graph and $\K$ any field. Then the following are equivalent:

\begin{enumerate}[\upshape(1)]
\item  Every cycle in $E$ has an exit. 
\item  The Leavitt path algebra $L_\K(E)$ is a Zorn ring.
\item The abelian group $\mathbb Z$ acts freely on the talented monoid $T_E$. 
\end{enumerate} 
\end{theorem}



We are now in position to state the two most important (and most investigated!),  as-yet-unresolved statements  in the theory of Leavitt path algebras. Recall that a directed graph is called 
\emph{strongly connected} if given any two distinct vertices $u$ and $v$,  there is  a directed path $p$, with $s(p)=u$ and $r(p)=v$. 

\begin{question}[Open since 2008] ``The Algebraic Kirchberg-Phillips Question"\label{conj111}
Let $E$ and $F$ be strongly connected finite graphs and $\K$ a field. 
Easily (2) implies (1) in the following.  The Question is whether (1) implies (2).

\begin{enumerate}[\upshape(1)]

\item There is a monoid isomorphism $\phi: M_{E} \rightarrow  M_{F} $, such that $\phi\big (\sum_{v\in E^0} v\big)=\sum_{v\in F^0} v$;

\smallskip

\item There is a $\K$-algebra  isomorphism $\varphi:L_\K(E) \rightarrow L_\K(F)$.
\end{enumerate}
\end{question}

\begin{conjecture}\label{HazratConjecture}(Open since 2013) \label{conj1}
Let $E$ and $F$ be finite graphs and $\K$ a field. 
Easily (2) implies (1) in the following. The Conjecture is that (1) implies (2).  

\begin{enumerate}[\upshape(1)]

\item There is a $\mathbb Z$-isomorphism $\phi: T_{E} \rightarrow  T_{F} $, such that $\phi\big (\sum_{v\in E^0} v\big)=\sum_{v\in F^0} v$;

\smallskip

\item There is a graded $\K$-algebra  isomorphism $\varphi:L_\K(E) \rightarrow L_\K(F)$.
\end{enumerate}
\end{conjecture}  

In fact both Question~\ref{conj111} and Conjecture~\ref{conj1} can be formulated in  terms of Grothendieck groups and thus  are related to the deep machinery of $K$-theory. One can show that the group completions of $M_E$  and $T_E$ are $K_0(L_\K(E))$ and $K^{\gr}_0(L_\K(E))$, the Grothendieck group and the graded Grothendieck group of $L_\K(E)$, respectively. 

A large portion of the 2024 Oberwolfach gathering ``Combinatorial $*$-algebras"  was focused on ideas surrounding Question~\ref{conj111} and Conjecture~\ref{conj1} \cite{oberwillie}.  For a  comprehensive account of the current status of these investigations, we refer the reader to~\cite{willie}.

\begin{remark}
There exist many generalizations of the notion of a Leavitt path algebra.  These include Leavitt path algebras of separated graphs~\cite{arasep}, weighted graphs~\cite{Preweight}, hypergraphs~\cite{Phyper},  higher rank graphs~\cite{kpalgebra},  Steinberg algebras~\cites{Steinberg1,clarkhaz}, and many more.    Each of these generalizations has had connections and applications to various ideas throughout mathematics. 
\end{remark}

\begin{remark}\label{L2motherRemark}
Referring to Question \ref{L2mother} above, it was shown in \cite{BS} that for every finite graph $E$ and every field $\K$, the Leavitt path algebra $L_\K(E)$ embeds in $L_\K(1,2)$.  

\end{remark}

\section{Connections}

\subsection{Connection with chip firing}

In \cite{AH}  the two authors discovered an unexpected connection between Leavitt path algebras and the previously discussed theory of chip firing \S\ref{chipsfire}.  All the ideas appearing in the statement of the following theorem  have already been presented in the current article.    

\begin{theorem} Let $E$ be a sandpile graph,  $\SP(E)$ its sandpile monoid, and $\K$ any field. Then $$\SP(E) \cong  \mathcal V(L_\K(E')),$$ 
where $E'$ denotes the vertex-weighted graph obtained from $E$.  
\end{theorem}

Although it looks quite cryptic, the theorem above relates chip firing with Leavitt path algebras. This opens the door for further understanding of the chip firing process by associating it to a very rich structure of algebras.  This is a demonstration of how different areas of mathematics can be related to each other and perhaps ultimately shed light on real life phenomena. 

\bigskip

\subsection{Connection with $C^*$-algebras}

In the first half of the 1930s in Berlin,  and in the second half of that same decade in Princeton, John von Neumann developed the theory of {\it algebras of operators}, in  part as a mathematical vehicle for quantum theory~\cite{vonneumannquant}. The theory of operator algebras has become a major branch of mathematics,  interacting with numerous other disciplines in mathematics, as well as several areas of theoretical physics. 
The $C^*$-algebras (self-adjoint operator algebras on Hilbert space which are closed in the norm topology) have been shaped by challenging problems that arise from physics, representation theory, and dynamical systems, and their study has  in return contributed to these subjects. 

The trend of algebraization of concepts from operator theory into a purely algebraic context started with von Neumann, Kaplansky and Jacobson,  who devised ways of seeing operator algebraic properties in underlying discrete structures.  As Berberian puts it in \cite{berber}, ``if all the functional analysis is stripped away ... what remains should stand firmly as a substantial piece of algebra, completely accessible through algebraic avenues''.

 Leavitt path algebras stand as a sparkling example of Berberian's point of view.    They are the discrete version of graph $C^*$-algebras, and are rich enough in structure to preserve a significant amount of  data, even without the presence of the topology (in particular the norm)  of graph $C^*$-algebras.  As an example, in the foundational work \cite{amp},  Ara, Moreno and Pardo show that $K_0(L_\C(E)) \cong K_0(C^*(E))$.  As another  example, the graph $C^*$-algebra $C^*(E)$ is simple (i.e., has no proper closed two-sided ideals)  if and only if the simplicity criteria for $L_\K(E)$ given in Theorem~\ref{firsttheorem} hold for the graph $E$. 
 
 Here is a fundamental  conjecture (the so called \emph{Abrams-Tomforde Conjecture},  posed in 2011  \cites{abramstomforde, eilers3}) that is germane to the current discussion.    It  relates the collection of Leavitt path algebras to the collection of graph $C^*$-algebras.
\begin{conjecture}\label{ATConjecture}  The following statements are equivalent for any graphs $E$ and $F$,  and any field $\K$.  

(1)  The Leavitt path algebras $L_\K(E)$ and $ L_\K(F)$  are isomorphic as rings.

(2)   The graph $C^*$-algebras $C^*(E)$ and $C^*(F)$ are isomorphic as $C^*$-algebras .
\end{conjecture}

\noindent 
The implication $(1) \Rightarrow (2)$ of Conjecture \ref{ATConjecture} was verified in 2021 for  graphs having finitely many vertices.   This was one of many important results established  in the extremely powerful, extremely deep, nearly-100-page-long work of Eilers, Restorff, Ruiz and S\o rensen \cite{ERRS}.  

With Conjecture \ref{ATConjecture} as context, we note that Conjecture~\ref{conj1} can be extended to include a third  statement about  graph $C^*$-algebras:  

(3)  there exists  a $\mathbb T$-equivariant $\ast$-isomorphism $C^*(E) \cong C^*(F)$. (Here $\mathbb T$ is the circle group, which induces a gauge action on $C^*(E)$ and $C^*(F)$.)

\subsection{Connection with group theory}
In 1965, in handwritten notes, Richard Thompson discovered two groups which are infinite and simple (i.e., have no nontrivial normal subgroups), and are also finitely presented (i.e., they are generated by a finite number of generators which satisfy  a finite number of relations).  Subsequently, in 1974 Higman constructed an infinite family of such groups:   specifically, for each pair of positive integers $(n,r)$, there is an infinite finitely presented simple group (denoted  $G_{n,r}^+$, and called the {\it Higman-Thompson group of type} $(n,r)$).  One of Thompson's two groups arises as  $G_{2,1}^+$.   Sixty years on, these groups are still the topics of intense research.

The following group-theoretic problem was open, and  intensely investigated,  for more than thirty years:     give  conditions on pairs of positive integers $(n,r)$ and $(m,s)$  which are equivalent to the existence of an isomorphism between the groups $G_{n,r}^+$ and $G_{m,s}^+$.   In 2011 Enrique Pardo showed in \cite{pardo}  that the group $G_{n,r}^+$ could be realized as an explicitly described set of  (two-sided) invertible elements inside the matrix ring ${\rm M}_r(L_\K(1,n))$.   Three years earlier,   Abrams, \'{A}nh and Pardo had solved the isomorphism question for matrix rings over Leavitt algebras, to wit:

{\bf Theorem} \cite[Theorem 5.2]{AAP}:   For pairs of positive integers $(n,r)$ and $(m,s)$, and any field $\K$, 
$${\rm M}_r(L_\K(1,n)) \cong {\rm M}_s(L_\K(1,m)) \  \ \Longleftrightarrow \ \ m=n \ \mbox{and} \ {\rm gcd}(r,n-1) = {\rm gcd}(s,n-1).$$
\noindent
(We note that the equation ${\rm gcd}(r,n-1) = {\rm gcd}(s,n-1)$ is equivalent to the condition that appears in \cite[Theorem 5.2]{AAP}.) 

Using the explicitly constructed  isomorphisms presented in \cite{AAP}, Pardo solved the longstanding  isomorphism question for the Higman-Thompson groups, to wit:  

{\bf Theorem} \cite[Theorem 3.6]{pardo}:    For pairs of positive integers $(n,r)$ and $(m,s)$, 
$$G_{n,r}^+ \cong G_{m,s}^+  \ \ \Longleftrightarrow \ \  m=n \ \mbox{and} \ {\rm gcd}(r,n-1) = {\rm gcd}(s,n-1).$$

\section{Acknowledgements}
The idea for this paper formed during the Combinatorial $*$-Algebras workshop at Oberwolfach in April 2024. The authors would like to thank Mathematisches Forschungsinstitut Oberwolfach for hosting the gathering. 
Hazrat acknowledges Australian Research Council Discovery Project DP230103184.

\end{document}